\documentclass{agtart_a}
\pdfoutput=1

\usepackage{dbnsymb}
\usepackage[all]{xy}
\usepackage{graphicx}
\usepackage{pinlabel}
\usepackage{color}

\title[The Karoubi envelope and Lee's degeneration of Khovanov homology]{The Karoubi envelope\\ and Lee's degeneration of Khovanov homology}

\author{Dror Bar-Natan}
\givenname{Dror}
\surname{Bar-Natan} 
\address{Department of Mathematics\\
University of Toronto\\\newline
Toronto Ontario M5S 2E4\\Canada}
\email{drorbn@math.toronto.edu}
\urladdr{http://www.math.toronto.edu/~drorbn}

\author{Kim Morrison}
\givenname{Kim}
\surname{Morrison}
\address{Department of Mathematics\\
University of California, Berkeley\\\newline
Berkeley CA 94720\\USA}
\email{scott@math.berkeley.edu}
\urladdr{http://math.berkeley.edu/~scott}

\volumenumber{6}
\issuenumber{}
\publicationyear{2006}
\papernumber{53}
\startpage{1459}
\endpage{1469}

\doi{}
\MR{}
\Zbl{}


\received{29 June 2006}
\revised{}
\accepted{20 July 2006}
\published{4 October 2006}
\publishedonline{4 October 2006}
\proposed{}
\seconded{}
\corresponding{}
\editor{}
\version{}

\arxivreference{math.GT/0606542}

\let\xysavmatrix\xymatrix
\def\xymatrix{\disablesubscriptcorrection\xysavmatrix}

\SelectTips{cm}{}

\makeatletter
\def\cnewtheorem#1[#2]#3{\newtheorem{#1}{#3}[section]
\expandafter\let\csname c@#1\endcsname\c@prop}
\makeatother

\theoremstyle{plain}
\newtheorem{prop}{Proposition}[section]
\cnewtheorem{thm}[prop]{Theorem}
\cnewtheorem{prob}[prop]{Problem}
\cnewtheorem{lem}[prop]{Lemma}
\theoremstyle{definition}
\cnewtheorem{defn}[prop]{Definition}         

\numberwithin{equation}{section}

\def\calC{{\mathcal C}}

\def\calI{{\mathcal I}}
\def\calO{{\mathcal O}}

\def\llbracket{\left[\mskip -1.4mu\left[}
\def\rrbracket{\right]\mskip -1.4mu\right]}
\newcommand{\Kh}[1]{\llbracket#1\rrbracket_0}
\newcommand{\KhL}[1]{\llbracket#1\rrbracket_1}

\newcommand{\Cobz}{{\mathcal Cob}_0}
\newcommand{\Cobo}{{\mathcal Cob}_1}

\newcommand{\Cobdl}{{\mathcal Cob}_{\bullet/l}}
\newcommand{\Cobl}{{\mathcal Cob}_{/l}}
\newcommand{\im}{\operatorname{im}}
\newcommand{\Kar}{\operatorname{Kar}}
\newcommand{\Kom}{\operatorname{Kom}}

\newcommand{\mor}{\operatorname{mor}}
\newcommand{\Mat}{\operatorname{Mat}}

\newcommand{\Integer}{\mathbb Z}

\newcommand{\Iso}{\cong}

\newcommand{\htpy}{\simeq}

\newcommand{\tensor}{\otimes}
\newcommand{\Tensor}{\bigotimes}
\newcommand{\Hom}[2]{\operatorname{Hom}\left(#1,#2\right)}

\newcommand{\eps}[2]{{\hspace{-3pt}\begin{array}{c}%
  \raisebox{-2.5pt}{\includegraphics[width=#1\textwidth]{\figdir/#2}}%
\end{array}\hspace{-3pt}}}

\begin{document}

\begin{asciiabstract}
We give a simple proof of Lee's result from [Adv. Math. 179 (2005)
554-586], that the dimension of the Lee variant of the Khovanov
homology of a c-component link is 2^c, regardless of the number of
crossings.  Our method of proof is entirely local and hence we can
state a Lee-type theorem for tangles as well as for knots and
links. Our main tool is the ``Karoubi envelope of the cobordism
category'', a certain enlargement of the cobordism category which is
mild enough so that no information is lost yet strong enough to allow
for some simplifications that are otherwise unavailable.
\end{asciiabstract}

\begin{webabstract}
We give a simple proof of Lee's result from [Adv. Math. 179 (2005) 
554--586], that
the dimension of the Lee variant of the Khovanov homology of a
$c$--component link is $2^c$, regardless of the number of crossings.
Our method of proof is entirely local and hence we can state a
Lee-type theorem for tangles as well as for knots and links. Our
main tool is the ``Karoubi envelope of the cobordism category'', a
certain enlargement of the cobordism category which is mild enough
so that no information is lost yet strong enough to allow for some
simplifications that are otherwise unavailable.
\end{webabstract}

\begin{htmlabstract}
We give a simple proof of Lee's result from [Adv. Math. 179 (2005) 554--586], that
the dimension of the Lee variant of the Khovanov homology of a
c&ndash;component link is 2<sup>c</sup>, regardless of the number of crossings.
Our method of proof is entirely local and hence we can state a
Lee-type theorem for tangles as well as for knots and links. Our
main tool is the &ldquo;Karoubi envelope of the cobordism category&rdquo;, a
certain enlargement of the cobordism category which is mild enough
so that no information is lost yet strong enough to allow for some
simplifications that are otherwise unavailable.
\end{htmlabstract}

\begin{abstract}
We give a simple proof of Lee's result from~\cite{MR2173845}, that
the dimension of the Lee variant of the Khovanov homology of a
$c$--component link is $2^c$, regardless of the number of crossings.
Our method of proof is entirely local and hence we can state a
Lee-type theorem for tangles as well as for knots and links. Our
main tool is the ``Karoubi envelope of the cobordism category'', a
certain enlargement of the cobordism category which is mild enough
so that no information is lost yet strong enough to allow for some
simplifications that are otherwise unavailable.
\end{abstract}

\maketitle

\section{Introduction} \label{sec:Introduction}
In a beautiful article~\cite{MR2173845}, Eun Soo Lee introduced a
second differential $\Phi$ on the Khovanov complex of a knot (or
link) and showed that the resulting double complex has
uninteresting homology. In a seemingly contradictory manner, this
is a very interesting result --- for this ``degeneration'' of the
Lee theory is in itself an extra bit of information about the
original Khovanov homology, masterfully used by
Rasmussen~\cite{math.GT/0402131} to define the aptly named
``Rasmussen invariant'' of a knot and to give a combinatorial proof
of an old conjecture of Milnor.

Unfortunately Lee's proof of her degeneration result is a bit
technical and inductive in nature.
The purpose of this note is to reprove Lee's degeneration result in
local terms, using tools in the spirit
of~\cite{Bar-Natan:Cobordisms}. Thus in addition to being a bit more
conceptual, our methods work for tangles as well as for knots and
links.

Another proof of Lee's result was found by S\,Wehrli
\cite{math.GT/0409328}. His proof is quick, short and elegant; the
only small advantage of our proof over his is its locality and
applicability to tangles.

Let us sketch our results now; the relevant terminology (which
closely follows~\cite{Bar-Natan:Cobordisms}) will be quickly
recalled in \fullref{sec:QuickReview} below.

A \emph{confluence} within a smoothing $S$ of a tangle $T$ is a pair
of arc segments in $S$ that correspond to a small neighbourhood of a
crossing in $T$. Thus if $T$ has $n$ crossings, $S$ will have $n$
confluences; a crossing such as $\slashoverback$ in $T$ becomes a
confluence such as $\smoothing$ or $\hsmoothing$ in $S$.

\begin{defn} Let $T$ be a tangle and let $S$ be a smoothing of $T$.
An alternate colouring of $S$ is a colouring of the components of
$S$ with two colours (always ``red'' and ``green'' below), so that
the two arc segments at every confluence of $S$ are coloured with
different colours.
\end{defn}

Our main theorem is the following:

\begin{thm} \label{thm:main}
Within the appropriate category (see below), the Khovanov--Lee
complex of a tangle $T$ is homotopy equivalent to a complex with one
generator for each alternately coloured smoothing of $T$ and with
vanishing differentials.
\end{thm}

Alternately coloured smoothings are easy to manage:

\begin{prop} \label{prop:main}
Any $c$--component tangle $T$ has exactly $2^c$ alternately coloured
smoothings. These smoothings are in a bijective correspondence with
the $2^c$ possible orientations of the $c$ components of $T$.
\end{prop}

Together \fullref{thm:main} and \fullref{prop:main}
imply Lee's result that the Khovanov--Lee homology of an
$c$--component link is $2^c$ dimensional.

What is that ``appropriate category'' of \fullref{thm:main}? It
is the category of complexes over the ``Karoubi envelope''
$\Kar(\Cobo)$ of the category $\Cobo$ used to describe the
Khovanov--Lee complex. It is perhaps the nicest gadget appearing in
our note --- it introduces more objects into $\Cobo$ allowing for
more opportunities to simplify complexes over $\Cobo$, yet it does
not introduce new morphisms between existing objects of $\Cobo$, and
hence no new homotopies or homotopy equivalences. Thus there is no
loss of information in the passage from $\Cobo$ to $\Kar(\Cobo)$.

\subsection{Is there anything left to do?} Plenty. Lee's
degeneration is critical to understanding
Rasmussen~\cite{math.GT/0402131}. We have a clean approach to Lee's
degeneration, and one may hope it will help improve our
understanding of Rasmussen's work. This still remains to be done.
More specifically, the Khovanov--Lee homology is filtered and the
definition of the Rasmussen invariant requires this filtration.
Thus:

\begin{prob} Figure out how the filtration in the Khovanov--Lee
homology interacts with everything down below.
\end{prob}

\subsection{The plan} In \fullref{sec:QuickReview} we quickly
recall the relevant definitions of the Khovanov and Khovanov--Lee
theories. Then in \fullref{sec:envelope} we review the definition of
the Karoubi envelope of a general category $\calC$.
\fullref{sec:RedAndGreen} is the heart of the paper. In it we
introduce the ``red'' and the ``green'' projections, the red-green
splitting within $\Kar(\Cobo)$ of objects in $\Cobo$ and prove
\fullref{thm:main}. Finally, in \fullref{sec:orientations} we
prove \fullref{prop:main}.

\subsubsection*{Acknowledgement} We wish to thank S\,Wehrli for drawing
our attention to~\cite{math.GT/0409328} and A\,Referee for further
comments and suggestions. This work was partially supported by NSERC
grant RGPIN 262178.

\section{A quick review of the local Khovanov theory}
\label{sec:QuickReview}

\begin{figure}[ht!]
\labellist
\small\hair 2pt
\pinlabel $00$ at 72 94
\pinlabel {$0^{\ast}$} [br] at 106 121
\pinlabel $01$ at 144 142
\pinlabel {$1^{\ast}$} [bl] at 179 122
\pinlabel $11$ at 215 94
\pinlabel {$1^{\ast}$} [tl] at 178 67
\pinlabel $10$ at 144 47
\pinlabel {$0^{\ast}$} [tr] at 107 68
\pinlabel $0$ at 72 3
\pinlabel $1$ at 144 3
\pinlabel $2$ at 215 3
\endlabellist
\[
  \begin{array}{c}
    \begin{picture}(0,0)%
\includegraphics{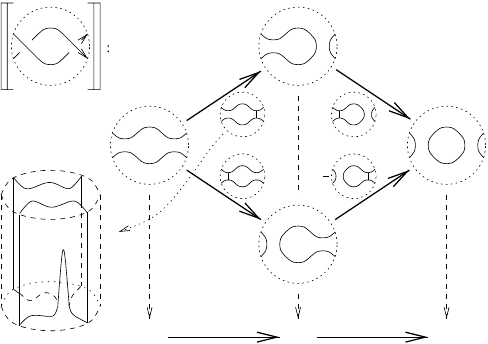}%
\end{picture}%
%
%
\setlength{\unitlength}{2605sp}%
\begingroup\makeatletter\ifx\SetFigFont\undefined%
\gdef\SetFigFont#1#2#3#4#5{%
  \reset@font\fontsize{#1}{#2pt}%
  \fontfamily{#3}\fontseries{#4}\fontshape{#5}%
  \selectfont}%
\fi\endgroup%
\begin{picture}(5894,4144)(1189,-3968)
\end{picture}%
  \end{array}
\]
\caption{The Khovanov complex of a $2$--crossing tangle}
\label{fig:KhovanovComplex}
\end{figure}
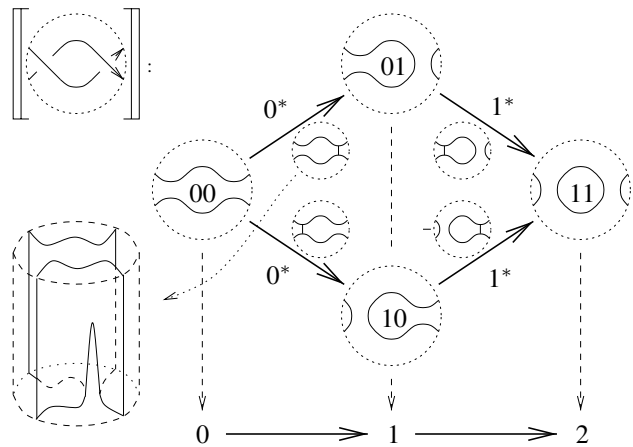

Let us briefly recall the definition of the Khovanov complex for
tangles, following \cite{Bar-Natan:Cobordisms}. Given an
$n$--crossing tangle $T$ with boundary $\partial T$ (such as the
$2$--crossing tangle in \fullref{fig:KhovanovComplex}) one
constructs an $n$--dimensional ``cube'' of $1$--dimensional smoothings
and $2$--dimensional cobordisms between them (as illustrated in
\fullref{fig:KhovanovComplex}). This cube is then ``flattened''
to a ``formal complex'' $\Kh{T}$ in the additive category
$\Cobz(\partial T)$ (denoted $\Cobdl^3(\partial T)$
in~\cite{Bar-Natan:FastKh,Bar-Natan:Cobordisms}) whose objects are
formally graded smoothings with boundary $\partial T$ and whose
morphisms are formal linear combinations of cobordisms whose tops
and bottoms are smoothings and whose side boundaries are
$I\times\partial T$, modulo some local relations. An overall height
shift ensures that the ``oriented smoothing'' appears in homological
height $0$.

The Khovanov complex $\Kh{T}$ is an object in the category
$\Kom(\Mat(\Cobz(\partial T)))$ of complexes of formal direct sums
of objects in $\Cobz(\partial T)$ and it is invariant up to
homotopies.

For simplicity we are using as the basis to our story one of the
simpler cobordism categories $\Cobz:=\Cobdl^3$ that appear
in~\cite{Bar-Natan:Cobordisms}, rather than the most general one,
$\Cobl^3$. It is worthwhile to repeat here the local relations that
appear in the definition of $\Cobz$
(see~\cite[Section 11.2]{Bar-Natan:Cobordisms}):
\begin{equation} \label{eq:LocalRelations}
\begin{array}{c}
  \begin{array}{c}
    \includegraphics[height=1cm]{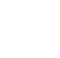}
  \end{array}\hspace{-2mm}=0,
  \qquad\qquad
  \begin{array}{c}
    \includegraphics[height=1cm]{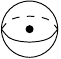}
  \end{array}\hspace{-2mm}=1,
  \qquad\qquad
  \begin{array}{c}\includegraphics[height=10mm]{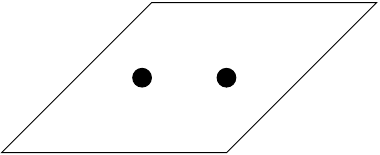}\end{array}
  \hspace{-4mm}=0,
\\
  \text{and}\qquad
  \begin{array}{c}\includegraphics[height=10mm]{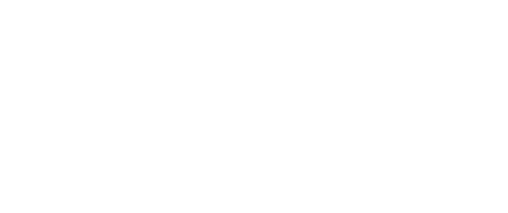}\end{array}
  =\begin{array}{c}\includegraphics[height=10mm]{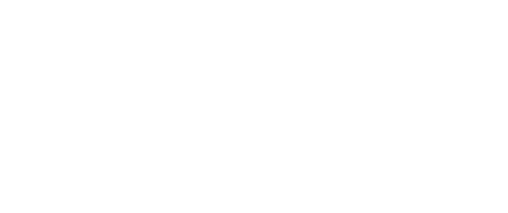}\end{array}
  +\begin{array}{c}\includegraphics[height=10mm]{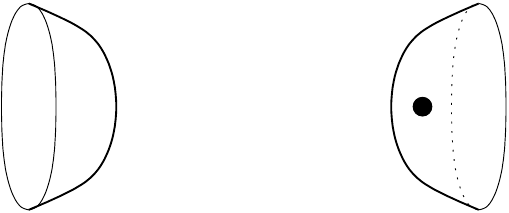}\end{array}.
\end{array}
\end{equation}
(If you're more used to the purely topological cobordism model
without dots, recall the translation $\eps{0.018}{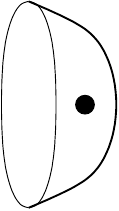} = \frac{1}{2}
\eps{0.04}{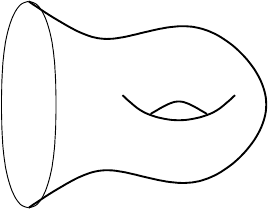}$. Note that while one of the original motivations
for describing the ``dotted'' theory was being able to work over
$\Integer$, we'll later need at least $2$ to be invertible and so
work over $\Integer_{(2)}$.)

Also recall from~\cite[Section 5]{Bar-Natan:Cobordisms} that
$\Kh{\,\cdot\,}$ is a planar algebra morphism. That is, if $T_1$ and
$T_2$ are tangles and $D(T_1,T_2)$ denotes one of their side-by-side
compositions (a side by side placement of $T_1$ and $T_2$ while
joining some of their ends in a certain way prescribed by a planar
arc diagram $D$), then $\Kh{D(T_1,T_2)}=\Kh{T_1} \tensor_D
\Kh{T_2}$. Here, as in~\cite[Section 5]{Bar-Natan:Cobordisms}, $\Kh{T_1}
\tensor_D \Kh{T_2}$ is the ``tensor product'' operation induced on
formal complexes by the horizontal composition operation $D$.

The Lee variant of Khovanov homology, or the Khovanov--Lee complex
$\KhL{T}$ of a tangle $T$, is constructed in exactly the same way as
$\Kh{T}$; the only difference is that the zero in the third local
relation in \eqref{eq:LocalRelations} is replaced with a one. Hence
it is valued in the category of complexes over the category
$\Cobo(\partial T)$ defined in exactly the same manner as
$\Cobz(\partial T)$, except with the following collection of local
relations (spot the one difference!):
\begin{equation} \label{eq:LeeRelations}
\begin{array}{c}
  \begin{array}{c}
    \includegraphics[height=1cm]{\figdir/S}
  \end{array}\hspace{-2mm}=0,
  \qquad\qquad
  \begin{array}{c}
    \includegraphics[height=1cm]{\figdir/Sd}
  \end{array}\hspace{-2mm}=1,
  \qquad\qquad
  \begin{array}{c}\includegraphics[height=10mm]{\figdir/ddot}\end{array}
  \hspace{-4mm}=1,
\\
  \text{and}\qquad
  \begin{array}{c}\includegraphics[height=10mm]{\figdir/CNN}\end{array}
  =\begin{array}{c}\includegraphics[height=10mm]{\figdir/dCNL}\end{array}
  +\begin{array}{c}\includegraphics[height=10mm]{\figdir/dCNR}\end{array}.
\end{array}
\end{equation}
While $\Cobz$ can be made into a graded category
(see~\cite[Section 6]{Bar-Natan:Cobordisms}), the degree of
$\begin{array}{c}\includegraphics[height=4mm]{\figdir/ddot}\end{array}$
is nonzero, so setting it to $1$ breaks the grading in $\Cobo$.
Otherwise the Khovanov theory and the Khovanov--Lee theories are
completely parallel. In particular, the Khovanov--Lee complex is also
an up-to-homotopy knot invariant and it is also a planar algebra
morphism.

\section{A quick review of the Karoubi envelope}
\label{sec:envelope}
A {\em projection\/} is an endomorphism $p$ satisfying $p^2=p$. In
many contexts in mathematics, if $p$ is a projection then so is
$1-p$, and together these two projections decompose space as a
direct sum of the image of $p$ with the image of $1-p$. Thus a
projection often gives rise to a decomposition of space into pieces
which are hopefully simpler.

The equation $p^2=p$ makes sense for an endomorphism in an arbitrary
category, so projections make sense in an arbitrary category. And if
$p$ is a projection then so is $1-p$ in an arbitrary {\em additive\/}
category. But in a general (additive or not) category, ``the image
of $p$'' (or of $1-p$) may or may not make sense.

The {\em Karoubi envelope\/}\footnote{The Karoubi envelope
construction \cite{wiki:Karoubi-Envelope} was first described by Freyd
\cite{MR0166240}, a few years before Karoubi. It has previously been
used in motivic cohomology by Mazur \cite{MR2104916} and in diagrammatic
representation theory, eg Kuperberg \cite{MR1403861}.} of a category $\calC$
is a way of adding objects and morphisms to $\calC$ so that every
projection has an image and so that if $p\co \calO\to\calO$ is a
projection and $\calC$ is additive, then (with the proper
interpretation) $\calO \Iso \im p\oplus\im (1-p)$. Thus sometimes
complicated objects can be simplified in the Karoubi envelope of
$\calC$, while in $\calC$ they may be indecomposable.

Let us turn to the formal definitions.

\begin{defn} Let $\calC$ be a category. An endomorphism
$p\co \calO\to\calO$ of some object $\calO$ in $\calC$ is called {\em a
projection} if $p\circ p=p$. The Karoubi envelope $\Kar(\calC)$ of
$\calC$ is the category whose objects are ordered pairs $(\calO, p)$
where $\calO$ is an object in $\calC$ and $p\co \calO\to\calO$ is a
projection. If $(\calO_1, p_1)$ and $(\calO_2, p_2)$ are two such
pairs, the set of morphisms in $\Kar(\calC)$ from $(\calO_1, p_1)$
to $(\calO_2, p_2)$ is the collection of all $f\co \calO_1\to\calO_2$
in $\calC$ for which $f=f\circ p_1=p_2\circ f$. An object $(\calO,
p)$ in $\Kar(\calC)$ may also be denoted by $\im p$.\footnote{If
you're worried about just introducing $\im p$ as notation, when you
already know a category-theoretic definition of image, eg
\cite{wiki:Image}, don't be; this \emph{is\/} actually an image.}
\end{defn}

The composition of morphisms in $\Kar(\calC)$ is defined in the
obvious way (by composing the corresponding $f$'s). The identity
automorphism of an object $(\calO, p)$ in $\Kar(\calC)$ is $p$
itself. It is routine to verify that $\Kar(\calC)$ is indeed a
category. There is an obvious embedding functor
$\calI\co \calO\mapsto(\calO,I)$ of $\calC$ into $\Kar(\calC)$ and
quite clearly, $\mor_{\Kar(\calC)}(\calI\calO_1,
\calI\calO_2)=\mor_\calC(\calO_1,\calO_2)$ for any pair of objects
$\calO_{1,2}$ in $\calC$. Thus we will simply identify objects in
$\calC$ with their image via $\calI$ in $\Kar(\calC)$.

Below we will assume that $\calC$ is an additive category and that
direct sums of objects make sense in $\calC$. As
in~\cite{Bar-Natan:Cobordisms}, there is no loss of generality in
making these assumptions as formal sums of morphisms and formal
direct sums of objects may always be introduced.

\begin{prop} Let $p\co \calO\to\calO$ be an endomorphism in $\calC$.
\begin{enumerate}
  \item If $p$ is a projection then so is $1-p$.
  \item In this case, $\calO \Iso \im p\oplus\im (1-p)$ in
  $\Kar(\calC)$.
\end{enumerate}
\end{prop}

\begin{proof}
\begin{enumerate}
  \item $(1-p)^2=1-2p+p^2=1-2p+p=1-p$ (sorry for the damage to the
  rainforest).
  \item The isomorphism $\calO\to\im p\oplus\im (1-p)$ is given by
  the $1\times 2$ matrix $\begin{pmatrix}p & 1-p\end{pmatrix}$.
  Its inverse is the $2\times 1$ matrix $\left(\begin{smallmatrix}p \\
  1-p\end{smallmatrix}\right)$.\proved
\end{enumerate}
\end{proof}

Observe that if $p$ is a projection on $\calO$ and $p'$ is a
projection on $\calO'$, then the set $\Hom{(\calO,p)}{(\calO',p')}$ may be
naturally identified with $p' \Hom{\calO}{\calO'} p$. In fact, even
before taking the Karoubi envelope, $\Hom{\calO}{\calO'}$ can be
expressed as a direct sum of $4$ ``matrix entries'', each obtained by
precomposing with $p$ or $1-p$, and postcomposing with $p'$ or
$1-p'$.

In this paper we are mainly interested in complexes whose ``chain
spaces'' are objects in some category $\calC$ as above. The previous
proposition tells us that there may be some gain by switching to
working over $\Kar(\calC)$ as we may have new decompositions of old
objects. The proposition below tells us that there is no loss of
information in doing so.

\begin{prop} Let $\Omega_1$ and $\Omega_2$ be complexes in
$\Kom(\calC)$. If $\Omega_1$ and $\Omega_2$ are homotopy equivalent
as complexes in $\Kom(\Kar(\calC))$, they were already homotopy
equivalent as complexes in $\Kom(\calC)$.
\end{prop}

\begin{proof}
A homotopy equivalence $\Upsilon$ in $\Kom(\Kar(\calC))$ between
$\Omega_1$ and $\Omega_2$ is a certain slew of morphisms between
objects appearing in $\Omega_1$ and objects appearing in $\Omega_2$
(a chain morphism going one way, another going the other way and a
couple of homotopies). All those morphisms are between objects in
$\calC$ and as noted above, $\Kar(\calC)$ introduces no new
morphisms between objects in $\calC$. So $\Upsilon$ is really in
$\calC$.
\end{proof}

\section{Red and green in Lee's theory} \label{sec:RedAndGreen}
Viewed from our perspective, the key to Lee's theorem is the
presence in $\Cobo$ of two complementary projections, the red
projection $r$ and the green projection $g$, that can be composed
both vertically (in the ``category'' direction of $\Cobo$) and
horizontally (in the ``planar algebra'' direction;
see~\cite[Section 8.2]{Bar-Natan:Cobordisms}). Let us start with some
elementary school algebra that contains all the calculations we will
need regarding $r$ and $g$.

\begin{lem} \label{lem:ElementaryAlgebra} Let $b$ (for \verb"\bullet") be a variable
satisfying $b^2=1$ (find it in~\eqref{eq:LeeRelations}!), let $r$
(for red) be $(1+b)/2$ and let $g$ (for green) be $(1-b)/2$. Then
$r$ and $g$ are the eigenprojections associated to the involution
$b$. In particular,
\begin{enumerate}
  \item $r$ and $g$ are projections: $r^2=r$ and $g^2=g$.
  \item $r$ and $g$ are complementary: $r+g=1$.
  \item $r$ and $g$ are disjoint: $rg=0$.
  \item $r$ and $g$ are eigenprojections of $b$: $br=r$ and $bg=-g$. \qed
\end{enumerate}
\end{lem}

In $\Cobo$, one may place a $b$ (ie a bullet) anywhere on any
cobordism and thus, treating linear combinations in the obvious
manner, one may place an $r$ or a $g$ anywhere on any cobordism. In
particular, we may place them on any ``vertical curtain'', ie on
any connected component of the identity morphism from a smoothing to
itself. The resulting ``identities labeled $r$ and/or $g$'' are
projections by the above lemma and hence they represent objects in
the Karoubi envelope $\Kar(\Cobo)$ of $\Cobo$. Thus in $\Kar(\Cobo)$
we have, for example, the isomorphism
\[
  \includegraphics[height=16mm]{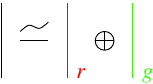}
\]
between a single arc smoothing and a direct sum of two single arc
smoothings, one paired with the $r$ projection and one with $g$.
Likewise we also have the isomorphism
\begin{equation} \label{eq:DecompositionExample2}
  \includegraphics[height=16mm]{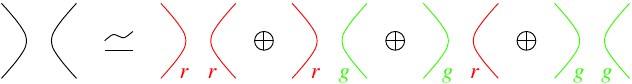}.
\end{equation}
We can now observe that a dramatic simplification occurs at the very
first step of calculating a knot invariant; the complex associated
to a single crossing now has an up-to-homotopy representative with a
vanishing differential. Indeed, $\KhL{\slashoverback}$ is the
two-step complex
$\xymatrix{\smoothing\ar[r]^{\smash{\HSaddleSymbol}}&\hsmoothing}$ in which
$\HSaddleSymbol$ denotes the saddle morphism. As
in~\eqref{eq:DecompositionExample2}, each of the two objects
($\smoothing$ and $\hsmoothing$) appearing in this complex becomes a
direct sum of four objects in $\Kar(\Cobo)$. The differential
$\HSaddleSymbol$ becomes a $4\times 4$ matrix $M$ all of whose
entries are saddles, and each such saddle carries $r$ and $g$
insertions to match the colourings of its domain and target
smoothings. But $r$ and $g$ are disjoint (see
\fullref{lem:ElementaryAlgebra}) and the saddle cobordism is
connected, so a saddle bearing insertions of more than one colour
vanishes and hence only two of the $16$ entries of $M$ survive. Thus
in $\Kom(\Kar(\Cobo))$,
\[
  \KhL{\slashoverback}
  \Iso\xymatrix@C=4cm{
    {\begin{bmatrix}
      {}_r\smoothing_r\\{}_r\smoothing_g\\{}_g\smoothing_r\\{}_g\smoothing_g
    \end{bmatrix}}
    \ar[r]^{\begin{pmatrix}
      \HSaddleSymbol_r & 0 & 0 & 0 \\
      0 & 0 & 0 & 0 \\
      0 & 0 & 0 & 0 \\
      0 & 0 & 0 & \HSaddleSymbol_g
    \end{pmatrix}}
    & {\begin{bmatrix}
      \hsmoothing^r_r \\ \hsmoothing^r_g \\ \hsmoothing^g_r \\ \hsmoothing^g_g
    \end{bmatrix}}
  }.
\]
This last complex is a direct sum of four complexes. The first and
the last of the four summands,
\[\xymatrix{{}_r\smoothing_r\ar[r]^\HSaddleSymbol_r&\hsmoothing_r^r}
\quad\text{and}\quad
\xymatrix{{}_g\smoothing_g\ar[r]^\HSaddleSymbol_g&\hsmoothing_g^g},\]
are contractible as $\HSaddleSymbol_r$ and $\HSaddleSymbol_g$ are
invertible (with inverses $\frac12\ISaddleSymbol_r$ and
$-\frac12\ISaddleSymbol_g$, respectively)\footnote{This follows from
the neck cutting relation (the last of~\eqref{eq:LeeRelations}) and
the final part of \fullref{lem:ElementaryAlgebra}.}. Thus up to
homotopy,
\[
  \KhL{\slashoverback}
  \htpy\xymatrix{
    {\begin{bmatrix}{}_r\smoothing_g\\{}_g\smoothing_r\end{bmatrix}}
    \ar[r]^0
    & {\begin{bmatrix}\hsmoothing^r_g \\ \hsmoothing^g_r\end{bmatrix}}
  },
\]
and as promised, we found a representative for
$\KhL{\slashoverback}$ with a vanishing differential.

\begin{proof}[Proof of \fullref{thm:main}]%
\fullref{thm:main} is now simply a matter of assembling the
pieces. The discussion above shows that it holds for tangles
consisting of a single crossing.

If we build a tangle $T$ by combining crossings $X_1$ through $X_n$
using a planar operation $D$, then $\KhL{T}$ is the tensor product
$\Tensor_D \KhL{X_i}$.

Since the complexes $\KhL{X_i}$ have an object for each of the four
alternately coloured smoothing of $X_i$ and no nonzero
differentials, the complex $\KhL{T}$ also has no differentials, but
at first sight too many objects. While every alternately coloured
smoothing of $T$ appears as an object (because an alternately
colouring smoothing can be divided into alternating coloured
smoothings of the constituent crossings), we also have smoothings
which are alternately coloured, but have different colours appearing
on arcs which are connected. Now remember that
\fullref{lem:ElementaryAlgebra} holds for horizontal compositions
as well as vertical ones. A strand with both colours inserted is the
same as a strand with the zero projection inserted and so is
equivalent in the category to the zero object; the extra objects all
disappear.
\end{proof}

\section{Alternately coloured smoothings and orienting components}
\label{sec:orientations}
In order to return a little closer to Lee's language
in~\cite{MR2173845}, we will now show that the alternately coloured
smoothings of \fullref{sec:RedAndGreen} are in a one-to-one
correspondence with orientations of the original tangle.

\begin{figure}[ht!]
\newcommand{\orientation}[1]{\eps{0.2}{bijection/whitehead_link_orientation#1}}
\newcommand{\sm}[1]{\eps{0.2}{bijection/whitehead_link_smoothing#1}}
\begin{align*}
 \orientation{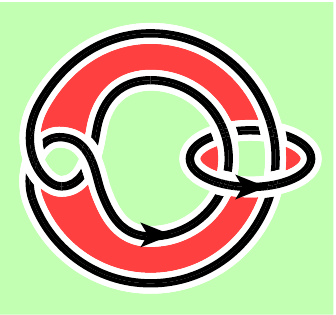} \leftrightarrow 
 & \labellist 
   \tiny\hair 2pt 
   \pinlabel {green} [l] at 0 6 
   \pinlabel {red} [r] at 160 10
   \endlabellist 
 \sm{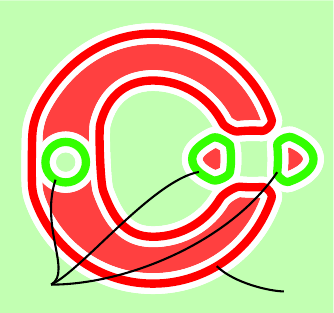}
 & \qquad
 \orientation{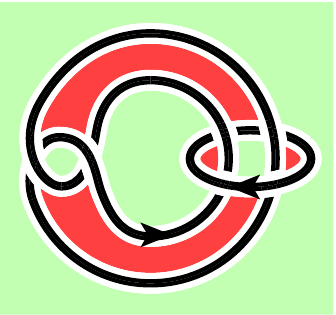} \leftrightarrow 
 & \labellist 
   \tiny\hair 2pt 
   \pinlabel {green} [l] at 0 6 
   \pinlabel {red} [r] at 161 10
   \endlabellist 
 \sm{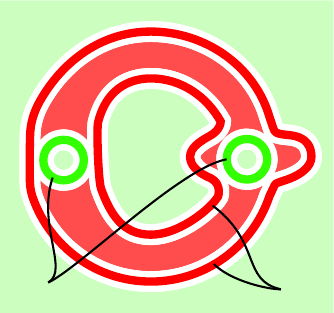} \\
 \orientation{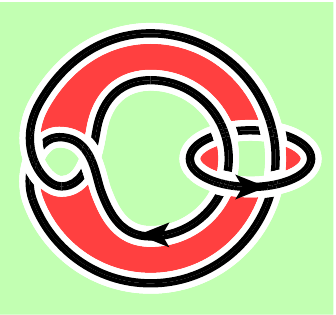} \leftrightarrow 
 & \labellist 
   \tiny\hair 2pt 
   \pinlabel {green} [l] at 0 6 
   \pinlabel {red} [r] at 154 10
   \endlabellist  
 \sm{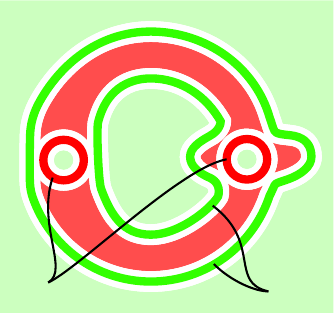} 
 & \qquad
 \orientation{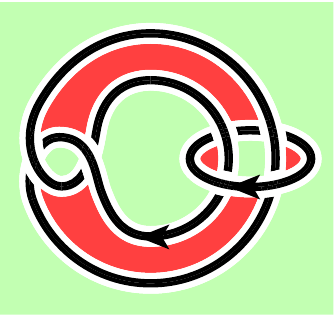} \leftrightarrow 
 & \labellist 
   \tiny\hair 2pt 
   \pinlabel {green} [l] at 0 5 
   \pinlabel {red} [r] at 150 10
   \endlabellist 
 \sm{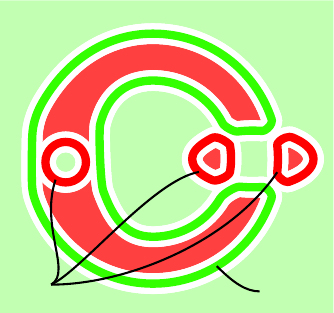}
\end{align*}
\caption{Generators of the Khovanov--Lee homology of the Whitehead
link. (In greyscale, replace ``red'' and ``green'' with ``dark'' and
``light''.) All four generators are in homological height zero,
because the two components have linking number $0$ (see
\fullref{prop:heights}).} \label{fig:whitehead_example}
\end{figure}

\begin{proof}[Proof of \fullref{prop:main}]
To begin, consider a tangle diagram $T$ and give its regions a
red-green checkerboard colouring. (The region outside the tangle
disk doesn't receive a colour; further, let's agree that the outer
region of a knot, or the marked boundary region of a tangle, is
green.)

We can canonically associate to an orientation of $T$ the oriented
smoothing, in which each arc is consistently oriented. We now need
to produce colours for the resulting arcs. As the arcs remain
oriented, we can simply take the colour appearing in the region on
the right. That this is consistent follows from the observation that
as you pass through a crossing, switching from one strand to the
other, the checkerboard colours appearing to your left and right
remain the same. Further, the colouring we've produced alternates
near each confluence, because as two incoming strands enter a
crossing they have opposite colours to their right.

Conversely, from an alternately coloured smoothing of a tangle, we
need to define an orientation of the tangle. Each arc of the knot
has been coloured either red or green, and it has a red and a green
region on either side. We orient the arc so that its own colour
appears on its right. Moving from one arc of the tangle to another
through a crossing, the colour of the arc changes, because the
smoothing is alternately coloured, and at the same time the
checkerboard colourings on either side switch. This ensures that two
opposite arcs at a crossing receive consistent orientations.

These two constructions are clearly inverses, so we have the desired
bijection. The four orientations of the Whitehead link and the
corresponding alternately coloured smoothings are shown in
\fullref{fig:whitehead_example}.
\end{proof}

Without any difficulty, we can also describe the homological height
of each generator. For this we need to compare a chosen orientation
of our tangle with the fixed ``original'' orientation (recall that the
orientation of a tangle is required in the definition of the
Khovanov homology to fix an overall homological height shift). Split
the tangle into two parts; $T_+$ containing those components where
the chosen and fixed orientations agree and $T_-$ containing the
components where the orientations disagree.

{
\newcommand{\lk}{\operatorname{lk}}
\begin{prop}
\label{prop:heights}%
Considering $T_+$ and $T_-$ to carry the
original orientations,
the homological height of the corresponding generator is
$\lk(T_+,T_-)$.
\end{prop}
\begin{proof}
This holds for the two component tangle consisting of a single
crossing. If both strands carry the original orientation, or both
carry the opposite orientation, then one of $T_+$ and $T_-$ is
empty, so the linking number is zero. The oriented smoothing of this
crossing is also the oriented smoothing with respect to the original
orientations and so sits in homological height zero. If exactly one
of the strands has been reversed, $\lk(T_+,T_-) = \pm 1$, agreeing
with the sign of the crossing (in the original orientation).
Happily, the oriented smoothing of the crossing is actually the
unoriented smoothing with respect to the original orientations, and
so sits in homological height $\pm 1$, again depending on the sign
of the crossing.

To extend the result to arbitrary tangles, it suffices to note that
every generator of $\KhL{T}$ is a planar composition of generators
for crossings and that both $\lk(T_+,T_-)$ and homological height
are additive under planar composition.
\end{proof}
}

\bibliographystyle{gtart}
\bibliography{link}

\end{document}